%% file: math0610037
\documentclass{birkmono}
\usepackage{bozhomac}	
\usepackage{bozhlogo}	
\usepackage{cite}	
\usepackage{shapepar}	
\usepackage{ifthen}	
\usepackage{array}	
\usepackage{supertab}	
\usepackage{tabularx}	
\usepackage{varioref}	
\usepackage{makeidx}	
\usepackage{index}	
\usepackage{pb-diagram}	
\usepackage{fancybox}   

\title{\bfseries	\vspace*{-2.2in}
\enlargethispage{0.55in}
{\huge Handbook of \\[1.1ex] normal frames and coordinates}\\[1.1ex]
}

\vspace{1.7ex}

\author{
Bozhidar Z. Iliev
\thanks{Laboratory of Mathematical Modeling in Physics,
Institute for Nuclear Research and \mbox{Nuclear} Energy,
Bulgarian Academy of Sciences,
Boul. Tzarigradsko chauss\'ee~72, 1784 Sofia, Bulgaria}
\thanks{E-mail address: bozho@inrne.bas.bg}
\thanks{URL: http://theo.inrne.bas.bg/$\sim$bozho/}
}

\date{
 \vspace{2.27ex}\ShortTitle{Handbook of normal frames and coordinates}
							\\[0.27ex]
\small
 \vspace{3.27ex}
	\begin{tabular}{r@{$\colon\to~$}l}
 \normalsize\sffamily\bfseries
Published in ``Progress in Mathematical Physics''&
 \normalsize\sffamily\bfseries
Vol. 42, Birkha\"user, Basel, 2006
\\
 \normalsize\sffamily\bfseries
&  Hardcover, ISBN: 3-7643-7618-X
\\[3ex]
 \normalsize\sffamily\bfseries
 \vspace{0.27ex} http://www.arXiv.org e-Print archive No. &
 \normalsize\sffamily\bfseries
 math.DG/0610037 \\[1.27ex]
 \vspace{0.27ex} Produced	& \fbox{\today}	\\[0.27ex]
	\end{tabular} \\[1.27ex]
	\begin{tabular}{r@{$\colon~$}l}
	\end{tabular} \\[-0.27ex]
  \vspace{4.27ex}{\Huge\BOZHO}	\\[4.27ex]
\vspace{0.27ex}
Subject Classes:
\Subject{Differential geometry}\\[2.27ex]
	\begin{tabular}{r@{\hspace{0.512em}}|@{\hspace{0.512em}}l}
 \vspace{0.27ex}
&
 \vspace{0.27ex}
	\end{tabular} \\[1.27ex]
 \vspace{0.27ex}\\[0.27ex]
}

 \makeindex		 	

\newcommand{\Label}[2][TWO*pi*Bozho*over*SEVEN]{
   \ifthenelse{\equal{TWO*pi*Bozho*over*SEVEN}{#1}}%
   {\label{#2}}%
   {\label{#1-#2}}}

\newcommand{\pageRef}[2][TWO*pi*Bozho*over*SEVEN]{
   \ifthenelse{\equal{TWO*pi*Bozho*over*SEVEN}{#1}}
   {\pageref{#2}}%
   {\pageref{#1-#2}}}

\newcommand{\ChapterNo}{0}

\providecommand{\Ref}[2][TWO*pi*Bozho*over*SEVEN]{}
 \renewcommand{\Ref}[2][TWO*pi*Bozho*over*SEVEN]{
   \ifthenelse{\equal{TWO*pi*Bozho*over*SEVEN}{#1}}
   {\ref{#2}}
   {\ifthenelse{\equal{Chapt#1}{Chapt\ChapterNo}}
	{\ref{#1-#2}}
	{\ref{Chapt#1}.\ref{#1-#2}}
   }%
}

\newcommand{\eRef}[2][TWO*pi*Bozho*over*SEVEN]{(\Ref[#1]{#2})}
\newcommand{\eReftag}[2][TWO*pi*Bozho*over*SEVEN]{
   \ifthenelse{\equal{TWO*pi*Bozho*over*SEVEN}{#1}}
   {\eRef{#2}}
   {\ifthenelse{\equal{Chapt#1}{Chapt\ChapterNo}}
	{\eRef[#1]{#2}}
	{\eref{#1-#2}\vpageref{#1-#2}}
   }%
}

\newcommand{\Reftag}[2][TWO*pi*Bozho*over*SEVEN]{
   \ifthenelse{\equal{TWO*pi*Bozho*over*SEVEN}{#1}}
   {\Ref{#2}}
   {\ifthenelse{\equal{Chapt#1}{Chapt\ChapterNo}}
	{\Ref[#1]{#2}}
	{\ref{#1-#2}\vpageref{#1-#2}}
   }%
}

\newcommand{\vRef}[2][TWO*pi*Bozho*over*SEVEN]{
   \ifthenelse{\equal{TWO*pi*Bozho*over*SEVEN}{#1}}
   {\vref{#2}}
   {\ifthenelse{\equal{Chapt#1}{Chapt\ChapterNo}}
	{\vref{#1-#2}}
	{\ref{Chapt#1}.\vref{#1-#2}}
   }%
}
\newcommand{\vpageRef}[2][TWO*pi*Bozho*over*SEVEN]{
   \ifthenelse{\equal{TWO*pi*Bozho*over*SEVEN}{#1}}
   {\vpageref{#2}}
   {\vpageref{#1-#2}}}

\newindex{symbols}{sdx}{snd}{Notation index\\[7ex]
\protect\parbox{\linewidth}{\normalsize\normalfont\indent
The frequently used symbols with a fixed meaning are listen below.
They are arranged in more or less alphabetical order. Since the sorting of
mathematical symbols is not unique, they are sorted according to up to three
criterions: the symbol's name or pronunciations (if any), the name of
symbol's kernel (root) letter, and the meaning of the symbol as a whole.
For some symbols such a classification is not unique, due to which they are
listen more than ones; \eg $\nabla$ can be found under the letters ``C'',
standing for Connection or Covariant derivative (sorting by meaning), and
``N'', standing for Nabla (sorting by name). Besides, a multiple appearance
of a symbol under different sorting letters may mean that it has different
meanings in different but similar contexts. The number(s) standing to the
right of the symbols denote the page(s) where the symbols first appear or
are defined.
\\[0.5ex]}}
 
\newindex{authors}{adx}{and}{Author index}
 
\newindex{subject}{idx}{ind}{Subject index\\[4ex]
\protect\parbox{\linewidth}{\normalsize\normalfont\indent
The \textbf{boldface} numbers refer to the page where a concept is
defined. The letter `n' after a page number means that the concept is
mentioned in a footnote. The ligature `ff' after a page number stands for
``and following pages''. Some important concepts are listen twice under
different names, \eg `\emph{coefficients} of linear connection' and
`linear connection \emph{coefficients of}'.
 \\
 } }

\renewcommand{\thechapter}{\Roman{chapter}}
\renewcommand{\thesection}{\arabic{section}}

\makeatletter
\@addtoreset{footnote}{section}
\makeatother

\begin{document}		

\include{nf-book0}

\chapter
[Manifolds, normal frames and Riemannian coordinates]
{Manifolds, normal frames\\ and Riemannian coordinates}
\markboth{\slshape CHAPTER~\thechapter.\hspace{0.75em}%
MANIFOLDS. RIEMANNIAN COORDINATES}{}

\renewcommand{\ChapterNo}{1}	

\section*{}
 \shapepar\nutshape{\bfseries\itshape
%
 	The basic differential-geometric concepts, such as differentiable
manifolds and mappings, tensors and tensor fields, and linear connections,
on which the book rests, are introduced. Partially the notation and
terminology employed are fixed. The normal frames and coordinates are defined
as ones in which the coefficients of a linear connection in them vanish on
some set. Certain their general properties are mentioned. The Riemannian
coordinates, which are normal at their origin, are described.
}\newpage

\section{Introduction}
\markright{\slshape \thesection.\hspace{0.75em}%
INTRODUCTION}

%
%
	The goal of this chapter is twofold: it introduces most of the basic
preliminary definitions and results on which our investigation rests
(sections~\Ref[1]{Sect2}, \Ref[1]{Sect3}, and~\Ref[1]{Sect4}) and it begins
the study of the normal frames and coordinates (sections~\Ref[1]{Sect3.4}
and~\Ref[1]{Sect5}).

	The main concepts of differential geometry required for the
understanding of the book are: differentiable manifolds and mappings,
submanifolds, Riemannian manifolds, tangent vectors and vector fields, tensors
and tensor fields, linear connections. The readers acquainted with them may only
look over the corresponding sections for our notation, omitting the major
text to which they may wish to return later, following the references to it.

	In more details, the contents of the chapter is as follows.

	The purpose of Sect.~\Ref[1]{Sect2} is to fix our terminology and
notation concerning differentiable manifolds and some typical to them natural
structures. This is not a summary of the differential geometry, only certain
basic concepts and particular relations between them required for our future
aims are presented. At first the concepts of topological and differentiable
manifolds are introduced, then tangent vectors, cotangent vectors, and tensors
and the corresponding fields of them on a manifold are defined. Also some
expressions in local bases (or frames) and coordinates are given. If the reader
is acquainted with all this, he/she  can simply look over this section for our
notation skipping the main text. A reader interested in deeper understanding
of these concepts, as well as in differential geometry as a whole, should
consult with the specialized literature. Here is a (random) selection of such
titles. An elementary introduction to differential geometry, with `physical'
orientation, can be found
in~\cite{Schutz,Isham,DNF-1,DNF-2,Gockeler&Schucker,Nash&Sen}. The same
purpose can serve the
books~\cite{Warner,Gromoll&et_al.,Brickell&Clark,Bishop&Crittenden} which are
more `mathematically' oriented. Our text follows the excellent
(text)books~\cite{K&N-1,Bruhat}. At last, the advanced
works~\cite{Sternberg,Lang/manifolds,Helgason,Greub&et_al.-1} can be
recommended. A brief synopsis of the mathematics preceding the introduction
of manifolds is given it~\cite{Brickell&Clark,Lang/manifolds,Greub&et_al.-1}
while~\cite{Bruhat,Teleman} contain an expanded presentation of the
`preliminary' to manifolds material. Of course, the reading of all of the
above\ndash mentioned serious books is not necessary for the understanding of
what follows. For this end, the reading of Sect.~\Ref[1]{Sect2} is sufficient and
the references cited may be consulted for more detains and proofs of some
assertions. The knowledge of the tensor analysis in coordinate\ndash dependent
language is desirable~\cite{Lovelock&Rund,Schouten/Ricci}. It is almost
sufficient for the most of this and subsequent chapters.

	In Sect.~\Ref[1]{Sect3}, we introduce the concept of linear connection
on a manifold. The approach chosen is, in a sense, middle between elementary
books on general relativity, such as~\cite{MTW,Weinberg}, and pure
mathematical ones on differential geometry,
like~\cite{K&N-1,Greub&et_al.-2}. We have tried to follow
closely~\cite{K&N-1,Bruhat,Brickell&Clark} but the abstracting material is
adapted to the goals of the present book. After a motivation for what the
connections are needed for, we introduce the linear connections via a system
of axioms for the covariant derivative of the algebra of tensor fields over a
given manifold. We employ this method since the theory of vector bundles, which
is not required for chapters~\ref{Chapt1}\Ndash\ref{Chapt3}, will be involve
into action only at the beginning of chapter~\ref{Chapt4}. In this connection,
let us mention that the linear connections can be defined only on the algebra
of vector fields on a manifold (\ie to the tangent to it bundle), and then they
admit a unique extension on the whole algebra of tensor fields~\cite[chapter~3,
proposition~7.5]{K&N-1}. A more advanced and deep treatment of the theory of
linear connections on manifolds and vector bundles can be found
in~\cite{Greub&et_al.-2,Bishop&Crittenden,Helgason,Sternberg,Poor}. We also
present the notion of a parallel transport (induced by a linear connection)
which will practically step on scene in chapter~\ref{Chapt4} but here is a
natural place for it to appear. It will be used in
chapters~\ref{Chapt1}\Ndash\ref{Chapt3} for proving and formulating some
results.  Sect.~\Ref[1]{Sect3} ends with a brief consideration of the
geodesics and exponential mapping.

	The concept of Riemannian metric and Riemannian connection are given
in Sect.~\Ref[1]{Sect4}. If the reader is interested in essence of Riemannian
geometry, he/she is referred, for example,
to~\cite{Schouten/Ricci,Eisenhart/Riemannian,Rashevskii,K&N-2,
Mishchenko&Fomenko,Gromoll&et_al.,Brickell&Clark,Bishop&Crittenden,
K&N-1,Bruhat}.

	In Sect~\Ref[1]{Sect3.4}, we introduce the main objects of our
investigation, the \emph{normal frames and coordinates}. We define them as
ones in which the coefficients of a linear connection vanish on a given set.
Some considerations on the uniqueness and (an)holonomicity of the normal
frames are presented too.

	Sect.~\Ref[1]{Sect5} contains a complete description of normal frames
at a given point on ($C^\infty$) Riemannian manifolds. This is done on the
base of Riemannian coordinates which turn to be normal at their
origin. The geodesic coordinates are pointed as other example of coordinates
normal at a point. Some general results, proved further in
chapter~\Ref{Chapt2}, concerning the existence of normal frames on
submanifolds are quoted. An expanded presentation of the problem of existence
of normal coordinates at a point of a $C^\infty$ Riemannian manifold is given
in~\cite{Eisenhart/Riemannian,Schouten/Ricci}, where also a list of original
early works on this topic can be found.

	In Sect.~\Ref[1]{Sect8} are presented a number of examples and
exercises of concrete Riemannian connection and coordinates/frames normal for
them on different sets. At first, the (locally) Euclidean and one\ndash
dimensional manifolds are considered. The (pseudo)spherical coordinates on
(pseudo)spheres are (partially) investigated for sets on which they are normal
for the Riemannian connection induced on them by the metric on them generated
by the Euclidean one of the Euclidean space in which the (pseudo)spheres are
embedded. Similar instance on the two dimensional torus is presented. The
cosmological models of Einstein, de~Sitter and Schwarzschild are considered (in
concrete coordinates) from the view\ndash point of normal frames/coordinates on
them. Some peculiarities of the light cone in Minkowski spacetime are pointed
too.

	Sect.~\Ref[1]{Sect7} deals with certain terminological problems
concerning bases and frames. Some links between these concepts are explicitly
formulated and/or derived.

	The chapter ends with some general remarks and conclusions in
Sect.~\Ref[1]{Conclusion}.


\include{nf-book1}

\chapter
[Existence, uniqueness and construction of\\ normal frames and coordinates
 for linear connections]
{Existence, uniqueness\\ and construction of normal frames and coordinates\\
 for linear connections}
\markboth{\slshape CHAPTER~\thechapter.\hspace{0.75em}%
NORMAL FRAMES FOR CONNECTIONS}{}


\renewcommand{\ChapterNo}{2}	

\addcontentsline{lot}{chapter}{Chapter~\thechapter}

\section*{}
\enlargethispage*{3ex} \vspace*{-7.5ex}
\heartpar{\bfseries\itshape
%
An in\ndash depth investigation of existence, uniqueness and construction of
frames and coordinates normal for linear connections on manifolds is given.
Detailed review of the literature dealing with normal coordinates is presented.
Some proofs are improved/generalized which entails a number of new results.
Similar problems in the case with non\ndash zero torsion are studied.
Main results: For arbitrary (resp.\ torsionless) connections frames
(resp.\ coordinates) normal at a single point and along path exist; they
exist on submanifolds of higher dimensions iff the parallel transport along
paths lying in them is path\ndash independent. Complete constructive
description of all, if any, frames and coordinates normal for arbitrary linear
connections.%
}
\newpage

\section{Introduction}
\markright{\slshape \thesection.\hspace{0.75em}%
INTRODUCTION}

	This chapter presents a complete exploration of the problems linked
to the existence, uniqueness, and construction of normal coordinates and
frames for manifolds endowed with a linear connection, with or without
torsion. The review of the literature dealing with normal coordinates is
mixed with new results. Such are, first of all, the ones concerning normal
frames, connections with non\ndash vanishing torsion, and the complete
constructive description of the  normal coordinates, if any.

	The methods for description of normal coordinates/frames on
Riemannian manifolds can \emph{mutatis mutandis} be transferred on arbitrary
manifolds, real or complex ($\field=\field[R],\field[C]$),~%
\footnote{~%
In the literature is often supposed $\field=\field[R]$ but this does not
influence the results.%
}
endowed with linear connection. The possibility for this is hidden in the
fact that the existence and properties of the normal coordinates/frames on a
Riemannian manifold is intrinsically connected with the properties of the
Christoffel symbols, \ie with the Riemannian connection, not with the
particular metric generating them. After this situation was clearly
understood, somewhere in
1922\Ndash1927~\cite{Veblen,Birkhoff,Veblen&Thomas-1923,Thomas-1925}
(see~\cite[p.~155]{Schouten/Ricci} for other references), the attention of
the mathematicians, working in the field, was completely switched to the
exploration of normal coordinates on manifolds with linear connections.
Practically only the symmetric (torsionless) case has been investigate (see
the comments after remark~\vRef[1]{Rem3.4}). Some random works,
like~\cite{Bortolotti,Pinl}, dealing with the asymmetric case
(non\ndash zero torsion) do not add nothing new as they simply note that the
symmetric parts~\eRef[1]{3.9} of the connection coefficients (in coordinate
frame) are coefficients of a symmetric linear connection to which the known
results for torsionless connections are applicable.

	Below in this chapter, in more or less modern terms and notation,
are reviewed all results concerning the existence of normal
coordinates/frames on manifolds endowed with symmetric linear connection. It
contains a number of original new results too.

	At first (Sect.~\Ref[2]{Sect6.1}), we concentrate on coordinates or
frames normal at a single point. We present the known classical methods in
this field~\cite{Schouten/Ricci,Schouten/physics,Lovelock&Rund} and then,
modifying the methods that will be given in chapter~\ref{Chapt3} in full
generality, we present a full description of these coordinates/frames.

	In Sect.~\Ref[2]{Sect6.2} the attention is turned on the coordinates
or frames normal along paths without self\ndash intersections. For symmetric
linear connections, we give a detailed description of the Fermi coordinates as
the first known coordinates of this kind with~\cite{Schouten/Ricci} being our
basic reference. Then, modifying the methods developed for similar but more
general problems (see chapter~\ref{Chapt3} and~\cite{bp-Frames-path}), we
derive a complete description of all coordinates or frames normal along
paths without self\ndash intersections or along locally injective paths in
manifolds with symmetric or, respectively, arbitrary linear connections.

	Several pages deal with problems concerning normal frames and
coordinates on submanifolds with maximum dimensionality
(Sect.~\Ref[2]{Sect6.3}), in particular on neighborhoods and on the whole
manifold. We prove that such frames or coordinates exist iff the connection
is (locally) respectively flat or flat and torsionless. A complete
description of the normal frames and coordinates in these cases is presented.
We also point to some links between normal frames and parallel transports for
flat linear connections.

	Section~\Ref[2]{Sect6.4} explores the problems of existence,
uniqueness, and construction of frames or coordinates normal on arbitrary
submanifolds. The classical results of~\cite{ORai} are reproduced in details
using modern notation.  Meanwhile, the corresponding proofs are improved,
some results are generalized for arbitrary connections, with or without
torsion, and new ones are presented. Next, we provide a complete constructive
description of all frames (resp.\ coordinates) normal on submanifolds of a
manifold with arbitrary (resp.\ torsionless) linear connection. Amongst a
number of general results, we prove that normal on a submanifold frames
(resp.\ coordinates) exist iff the parallel transport is path\ndash
independent along paths lying entirely in it (resp.\ and the connection is
torsionless).

	Section~\Ref[2]{Sect7} contains instances and exercises illustrating the
general theory of this chapter. Explicit expressions for frames and
coordinates normal at a single point in and along a great circle on a two
dimensional sphere are presented in a case of the Riemannian connection induced
from the Euclidean space in which the sphere is embedded. Some problems
connected with frames/coordinates normal for Weyl connections are investigated.
All frames/coordinates normal in the one\ndash dimensional case are explicitly
described. Similar problem is solved along a geodesic path in 2\ndash
dimensional manifold. All coordinates normal at a point in Einstein\ndash
de~Sitter spacetime are found.

	A brief recapitulation of the above items can be found in
Sect~\Ref[2]{Conclusion}.

\include{nf-book2}

\chapter
[\hspace*{0.1em}
Normal frames and coordinates\\\hspace*{0.1em}
			for derivations on differentiable manifolds]
{Normal frames and coordinates\\ for derivations\\ on differentiable manifolds}
\markboth{\slshape CHAPTER~\thechapter.\hspace{0.75em}%
NORMAL FRAMES ON MANIFOLDS}{}

\renewcommand{\ChapterNo}{3}	

\section*{}
\enlargethispage*{3ex} \vspace*{-6ex}
\shapepar\nutshape{\bfseries\itshape
	The existence, uniqueness, and construction of frames and coordinates
normal for derivations (along vector fields, fixed vector field, paths, and
fixed path) of the tensor algebra over a manifold are explored in details.
For arbitrary vector fields or paths, normal frames (resp.\
coordinates) exist always (resp.\ if the torsion vanishes); on other
submanifolds or along more general mappings necessary and sufficient
conditions for such existence  are derived. For derivations along fixed vector
field or path normal frames and coordinates exist always. With a few
exceptions, a complete constructive description of the normal frames and
coordinates, if any, is presented. Frames simultaneously normal for two
derivations are studied. With respect to the normal frames, the unique role
of the linear connections amongst the other derivations is pointed out.%
}\newpage

\section{Introduction}
\markright{\slshape \thesection.\hspace{0.75em}%
INTRODUCTION}

	The aim of this chapter is the investigation of frames and
coordinates normal for different kinds of derivations of the tensor algebra
over a differentiable manifold. Since the linear connections are a particular
example of such derivations, the presented here material is a direct
continuation and generalization of the one in chapter~\ref{Chapt2}. But, as
we shall see, a number of problems concerning normal frames and charts for
general derivations are `locally' reduced to the same problems for linear
connections and, consequently, their (local) solutions could be found, in
more or less ready form, in chapter~\ref{Chapt2}.

	Some of the results in the present chapter are partially based on
the ones in the series of
works~\cite{bp-Bases-n+point,bp-Bases-path,bp-Bases-general,
bp-Frames-n+point,bp-Frames-path,bp-Frames-general,bp-NF-D+EP}
and are completely revised and generalized their versions. But most of the
material is new and original.

	Sect.~\Ref[3]{Sect2} has an introductory character. The concepts of
derivations and derivations along vector fields of the tensor algebra over a
manifold are introduced. Their components, coefficients (if they are linear),
curvature, and torsion are defined. Next, in section~\Ref[3]{Sect3}, the normal
frames and coordinates are defined as ones in which the components of a
derivation along vector fields vanish (on some set). The equations describing
the transition to normal frames or coordinates are derived and the linearity
of a derivation along vector fields is pointed as a necessary conditions for
their existence.

	In Sect.~\Ref[3]{Sect4} (resp.\ Sect.~\Ref[3]{Sect5}) is proved that
at a single point (resp.\ along a (locally injective) path) frames normal for
a linear at it (resp.\ along it) derivation along vector fields always exist
and their complete descriptions are given. Besides, if the derivation is
torsionless, all normal coordinates are found. In
sections~\Ref[3]{Sect6}\Ndash\Ref[3]{Sect8}, the problems of existence,
uniqueness, and complete description of frames and local charts (or
coordinates) on neighborhoods, on submanifolds, and along (injective or
locally injective) mappings, respectively, for derivations along vector
fields are studied in details and solved.

	To the problems concerning frames or coordinates normal for derivations
along fixed vector field is devoted Sect.~\Ref[3]{Sect9}. The existence of
normal frames and coordinates in this case is proved. A complete description
of the frames normal at a single point, along a path, and on the whole
manifold are presented. The local charts (or coordinates) normal at a point
are completely described.  Along a path the explicit system of
differential equations, which the normal coordinates must satisfy and which
always have (local) solutions, is derived. A method for obtaining the
coordinates (locally) normal on the whole manifold is pointed in the $C^\infty$
case.

	Normal frames for derivations along paths are investigated in
section~\Ref[3]{Sect10}. After the introduction of the basic definitions and
notation, it is proved that frames normal for a derivation along a given
(fixed) path always exist and their general form is found. A (local)
holonomic extension of such frames, as well as of any other frame defined only
along a path, is constructed. For derivations along arbitrary paths is proved
that they admit normal frames iff they are covariant derivatives along paths
induced by linear connections for which normal frames exist. Since the normal
frames for the derivations and connections turn to be identical, all problems
for these frames are transferred to similar ones considered in
chapter~\ref{Chapt2}.

	Section~\Ref[3]{Sect11} deals with problems connected with frames
simultaneously normal for two derivations along arbitrary/fixed vector field
or path. Necessary and sufficient conditions for the existence of such frames
are found. In particular, in the case of arbitrary vector field or path, they
exist iff the two derivations coincide. Normal frames for mixed linear
connections are explored. It is shown that this range of problems is completely
and equivalently reduced to similar one for two, possibly identical, linear
connections, the contra- and co\ndash variant `parts' of the initial mixed
connection.

	In section~\Ref[3]{Sect12} are collected and commented some results
concerning linear connections obtained in the preceding sections of this
chapter.

	Section~\Ref[3]{Sect14} illustrates the theory of the preceding
sections with concrete examples.

	Section~\Ref[3]{Sect13} contains a discussion of some terminological
problems linked to the normal frames or coordinates.

	The chapter ends with certain general remarks in
Sect.~\Ref[3]{Conclusion}.


\include{nf-book3}

\chapter
[Normal frames in vector bundles]
{Normal frames\\ in vector bundles}
\markboth{\slshape CHAPTER~\thechapter.\hspace{0.75em}%
NORMAL FRAMES IN VECTOR BUNDLES}{}


\renewcommand{\ChapterNo}{4}	

\section*{}
\heartpar{\bfseries\itshape
	The theory of linear transports along paths in vector bundles,
generalizing the parallel transports generated by linear connections, is
developed. The normal frames for them are defined as ones in which their
matrices are the identity one. A number of results, including theorems of
existence and uniqueness, concerning normal frames are derived. Special
attention is paid to the case when the bundle's base is a manifold. The normal
frames are defined and investigated also for derivations along paths and along
tangent vector fields in the last case. Frames normal at a single point or
along a given path  always exist. On other subsets normal frames exist only in
the curvature free case. The privileged role of the parallel transports is
pointed out in this context.
}\newpage

\section {Introduction}
\markright{\slshape \thesection.\hspace{0.75em}%
INTRODUCTION}

	The analysis of corollary~\vRef[2]{Cor6.4.1} reveals that the
properties of the parallel transport assigned to a linear connection, not
directly the ones of the connection itself, are responsible for the existence
of frames normal on a submanifold for the connection.~%
\footnote{~%
Here the situation is similar to the one described in the second paragraph of
Sect.~\vRef[2]{Introduction}: the properties of the Christoffel symbols, not
directly the ones of the Riemannian metric generating them, are fully
responsible for the existence of coordinates normal at a single point in a
Riemannian manifold.%
}
This observation forms the groundwork of the idea the `normal' frames to be
defined directly for (parallel) transports without referring to the concept of
a (linear) connection (or some other derivation along vector fields). The main
obstacle for the realization of such an approach to `normal frames' is that,
ordinary, the concept of a parallel transport is a secondary one, it is
introduced on the base of the concept of a (linear) connection. To the
solution of the last problem and the development of the mentioned approach to
normal frames (in finite dimensional vector bundles) is devoted the present
chapter of the book. As we shall demonstrate below, the consistent realization
of the above idea leads to a completely new look on the `normal frames', which
is self\ndash contained and incorporates as special cases all of the results of
the preceding chapters.

	The material in sections~\Ref[4]{Sect2}--\Ref[4]{Sect5}
and~\Ref[4]{Sect6} is based on the work~\cite{bp-NF-LTP} and the one after
them is practically new and written especially for the present book.~%
\footnote{~%
Although, some initial ideas and results are borrowed from the
papers~\cite{bp-LTP-Cur+Tor,bp-LTP-Cur+Tor-prop,bp-LT-Deriv-tensors}.%
}

	In the present chapter is studied a wide range of problems
concerning frames normal for linear transports and derivations along paths in
vector bundles and for derivations along tangent vector fields in the case
when the bundle's base is a differentiable manifold. In the last case, when
tangent bundles are concerned, the only general result, known to the author and
regarding normal frames, is~\cite[p.~102, theorem~2.106]{Poor}.

	The structure of this chapter is as follows.

	Sect~\Ref[4]{Sect8} introduces some basic concepts from the theory of
(fibre) bundles, in particular of the one of vector bundles, required for the
investigations in this chapter. After the concepts of bundle, section, and
vector bundle are fixed, a special attention to the ones of liftings of paths
and derivations along paths, which will play an important role further, is
paid.  The tensor bundles over a manifold are pointed as particular examples
of vector bundles. Details on these and many other concepts regarding
(fibre) bundles, the reader can fined in the
monographs~\cite{Steenrod,Husemoller,Sze-Tsen,James,Greub&et_al.,K&N-1,Warner}.

	Sect.~\Ref[4]{Sect2} is devoted to the general theory of linear
transports along paths in vector fibre bundles which is a far reaching
generalization of the theory of parallel transports generated by linear
connections.~%
\footnote{~%
This section is based on the early
works~\cite{bp-LT-tensors-S,bp-LT-tensors-I,bp-LT-Deriv-tensors,
bp-LTP-general,bp-LTP-appl,bp-TP-general,bp-TP-parallelT} of the author. For
some more general results, see chapter~\Ref{Chapt5}.%
}
The general form and other properties of these transports are studied. A
bijective correspondence between them and derivations along paths is
established. In Sect.~\Ref[4]{Sect3}, the normal frames are defined as ones in
which the matrix of a linear transport along paths is the identity (unit) one
or, equivalently, in which its coefficients, as defined in
Sect.~\Ref[4]{Sect2}, vanish `locally'. A number of properties of the normal
frames are found. In Sect.~\Ref[4]{Sect4} is explored the problem of
existence of normal frames.  Several necessary and sufficient conditions for
such existence are proved and the explicit construction of normal frames, if
any, is presented.

	Sect.~\Ref[4]{Sect5} concentrates on, possibly, the most important
special case of frames normal for linear transports or derivations along
smooth paths in vector bundles with a differentiable manifold as a base.
A specific necessary and sufficient condition for existence of normal frames
in that case is proved. In particular, normal frames may exist only for those
linear transports or derivations along paths whose (2\Ndash index)
coefficients linearly depend on the vector tangent to the path along
which they act. Obviously, this is a generalization of the derivation
along curves assigned to a linear connection.  Sect.~\Ref[4]{Sect6} is devoted
to problems concerning frames normal for derivations along tangent vector
fields in a bundle with a manifold as a base. Necessary and sufficient
conditions for the existence of these frames are derived. The conclusion is
made that there is a one\nobreak-to\nobreak-one onto correspondence between
the sets of linear transports along paths, derivations along paths, and
derivations along tangent vector fields all of which admit normal frames.

	In the first part of Sect.~\Ref[4]{Sect7}, based
on~\cite{bp-LTP-Cur+Tor}, the concept of a curvature of a linear transport
along paths is introduced and some its properties are explored. In its second
part, relations between the curvature of a linear transports along paths and
the frames normal for them are studied. The main result is that only the
curvature free transports admit normal frames. The concept of a torsion of a
linear transport along paths in the tangent bundle over a manifold is
introduced in Sect.~\Ref[4]{Sect9} (cf. the early
paper~\cite{bp-LTP-Cur+Tor}). Links between the torsion and holonomic
normal frames are investigated. The vanishment of the torsion is pointed as a
necessary and sufficient condition for existence of normal coordinates on
submanifolds. If such coordinates exist, their complete description is given.

	Sect.~\Ref[4]{Sect10} deals with parallel transports in the tangent
bundles over manifolds and frames normal for these transports. It is shown
that the parallel transport assigned to a linear connection is a special kind
of a linear transport in tangent bundles. As a side result, an axiomatic
definition of a parallel transport is obtained, on the base of which a new
definition of a linear connection, equivalent to the usual one, is given. The
flat parallel transports are pointed as the only linear transports along paths
in tangent bundles which transports admit normal frames. The coordinates
normal for flat and torsionless parallel transports are explicitly presented.

	Sect.~\Ref[4]{Sect12} concerns a special type of normal frames in
which the 3\ndash index coefficients, if any, of a linear transport along
paths vanish.

	Sect.~\Ref[4]{Sect11} is similar to Sect.~\Ref[4]{Sect10}, but it
deals with the interrelations between different types of derivations along
vector fields over a manifold and the linear transports along paths in the
tangent bundle over it. As examples, particular derivations or transports, such
as Fermi\ndash Walker, Jaumann, etc., are considered.

	The aim of Sect.~\Ref[4]{13Sect} is twofold. On one hand
(Subsections~\Ref[4]{131Subsect1}\Ndash\Ref[4]{133Subsect3}), the rigorous
relations between  the theory of linear transports along paths in vector
bundles and the one of parallel transports and connections in these bundles
are investigated. On the base of the axiomatic approach to the theory of
parallel transports, as presented in~\cite{Poor}, we show how it (and hence
the one of connections) is incorporated as a special case in the general
theory of linear transports along paths. On another
hand (subsections~\Ref[4]{134Subsect4} and~\Ref[4]{135Subsect5}), we
demonstrate how the results concerning normal frames and derived for linear
connections on manifolds and linear transports along paths are almost
\emph{in extenso} applicable to the theory of parallel transports and
connections on vector bundles.

	In Sect.~\Ref[4]{14Sect} is introduced the notion of autoparallel
paths in manifolds whose tangent bundle is endowed with a linear transport
along paths. If this transport is a parallel one, it is proved that the
autoparallels coincide with the geodesics of the linear connection generating
the transport.

	The chapter ends with some notes in Sect.~\Ref[4]{Conclusion}.

\vspace{1ex}

	All fibre bundles in this chapter are vectorial ones. The
base and total bundle space of such bundles can be general topological
spaces. However, if some kind of differentiation in one/both of these spaces
is needed to be introduced (considered), it/they should possess a smooth
structure; if this is the case, we require it/they to be smooth, of class
$C^1$, differentiable manifold(s). Starting from Sect.~\Ref[4]{Sect5}, the
base and total bundle space are supposed to be $C^1$ manifolds.
Sections~\Ref[4]{Sect2}--\Ref[4]{Sect4} do not depend on the existence of a
smoothness structure in the bundle's base. Smoothness of the bundle space is
partially required in sections~\Ref[4]{Sect8}--\Ref[4]{Sect4}.~%
\footnote{~%
The bundle space is required to be a $C^1$ manifold in Sect~\Ref[4]{Sect8}
(starting from definition~\Ref[4]{Defn8.2}), in definition~\Ref[4]{Defn3.1'},
in proposition~\Ref[4]{Prop3.1}--\Ref[4]{Cor3.1**}, if~\eRef[4]{3.1c}
and~\eRef[4]{3.1d} are taken into account, in theorem~\Ref[4]{Thm3.2}, and in
proposition~\Ref[4]{Prop4.4}.%
}


\include{nf-book4}

\chapter
[Normal frames for connections on differentiable fibre bundles]
{Normal frames for connections\\ on differentiable fibre bundles}
\markboth{\slshape CHAPTER~\thechapter.\hspace{0.75em}%
NORMAL FRAMES FOR CONNECTIONS ON BUNDLES}{}

\renewcommand{\ChapterNo}{5}	

\addcontentsline{lot}{chapter}{Chapter~\thechapter}
\addcontentsline{lof}{chapter}{Chapter~\thechapter}

\section*{}
 \shapepar\nutshape{\bfseries\itshape
%
The general connection theory on differentialble fibre bundles, with emphasis
on the vector ones, is partially considered. The theory of frames normal for
general connections on these bundles is developed. Links with the theory of
frames normal for linear connections in vector bundles are revealed. Existence
of bundle coordinates normal at a given point and/or along injective horizontal
path is proved and a necessary and sufficient condition of existence of bundle
coordinates normal along injective horizontal mappings is proved. The concept
of a transport along paths in differentiable bundles is introduced. Different
links between connections, parallel transports (along paths) and transports
along paths are investigated.
}
\newpage

\section {Introduction}
\markright{\slshape \thesection.\hspace{0.75em}%
INTRODUCTION}

	All connections considered until now, on manifolds and on vector
bundles, were linear. It is well known that there exist non\ndash linear
connections on vector bundles as well as on non\ndash vector ones. Can normal
frames (and/or coordinates) be introduced for such more general connections? The
positive solution of that problem is the main goal of the present chapter of
this book. For the purpose and for a comparison with the definitions and
results already obtained is required some preliminary material on general
connection theory on differentiable  bundles, which is collected in
sections~\Ref[5]{Sect2}--\Ref[5]{Sect6}. On its base, the normal frames for
connections on such bundles are studied in sections~\Ref[5]{Sect7}
and~\Ref[5]{SectA}.

	Sections~\Ref[5]{Sect2}--\Ref[5]{Sect6} follow the
work~\cite{bp-ConTheo-review},~%
\footnote{~%
The presentation of the material in sections~\Ref[5]{Sect2}--\Ref[5]{Sect4}
is according to some of the main ideas of~\cite[chapters~1 and~2]{Rahula}, but
their realization here is quite different and follows the modern trends in
differential geometry.%
}
sections~\Ref[5]{Sect7} and~\Ref[5]{SectA} are a
slightly revised version of~\cite{bp-NF-GC}, and section~\Ref[5]{Sect8}
reproduces in a modified form the paper~\cite{bp-C-TP}

	The work is organized as follows.

	In Sect.~\Ref[5]{Sect2} is collected some introductory material, like
the notion of Lie derivatives and distributions on manifolds, needed for our
exposition. Here some of our notation is fixed too.

	Section~\Ref[5]{Sect3} is devoted to the general connection theory on
bundles whose base and bundles spaces are differentiable manifolds. From
different view\ndash points, this theory can be found in many works, like%
~\cite{Kobayashi-1957,K&N-1,Sachs&Wu,Nash&Sen,Warner,Bishop&Crittenden,
Yano&Kon,Steenrod,Sulanke&Wintgen,Bruhat,Husemoller,Mishchenko,R_Hermann-1,
Greub&et_al.-1,Atiyah,Dandoloff&Zakrzewski,Tamura,Hicks,Sternberg,Rahula,
Mangiarotti&Sardanashvily}.
In
Subsect.~\Ref[5]{Subsect3.1} are reviewed some coordinates and frames/bases on
the bundle space which are compatible with the fibre structure of a bundle.
Subsect.~\Ref[5]{Subsect3.2} deals with the general connection theory. A
connection on a bundle is defined as a distribution on its bundle space which
is complimentary to the vertical distribution on it. The notions of parallel
transport generated by connection and of specialized frame are introduced.
The fibre coefficients and fibre components of the curvature of a connection
are defined via part of the components of the anholonomicity object of a
specialized frame. Frames adapted to local bundle coordinates are introduced
and the local (2\ndash index) coefficients in them of a connection are
defined; their transformation law is derived and it is proved that a
geometrical object with such transformation law uniquely defines a connection.

	In Sect.~\Ref[5]{Sect4}, the general connection theory from
Sect.~\Ref[5]{Sect3} is specified on vector bundles. The most important
structures in/on them are the ones that are consistent/compatible with the
vector space structure of their fibres. The vertical lifts of sections of
a vector bundle and the horizontal lifts of  vector fields on its base are
investigated in more details in Subsect.~\Ref[5]{Subsect4.1}.
Subsect.~\Ref[5]{Subsect4.3} is devoted to linear connections on vector
bundles, \ie connections such that the assigned to them parallel transport is a
linear mapping. It is proved that the 2\ndash index coefficients of a linear
connection are linear in the fibre coordinates, which leads to the introduction
of the (3\ndash index) coefficients of the connection; the latter coefficients
being defined on the base space. The transformations of different objects under
changes of vector bundle coordinates are explored. The covariant derivatives
are introduced and investigated in Subsect.~\Ref[5]{Subsect4.4}. They are
defined via the Lie derivatives~\cite{Rahula} and a mapping realizing an
isomorphism between the vertical vector fields on the bundle space and the
sections of the bundle. The equivalence of that definition with the widespread
one, defining them as mappings on the module of sections of the bundle with
suitable properties, is proved. In Subsect.~\Ref[5]{Subsect4.5}, the affine
connections on vector bundles are considered briefly.

	In Section~\Ref[5]{Sect6}, some of the results of the previous sections
are generalized when frames more general than the ones generated by local
coordinates on the bundle space are employed. The most general such frames,
compatible with the fibre structure, and the frames adapted to them are
investigated. The main differential\ndash geometric objects, introduced in the
previous sections, are considered in such general frames. Particular attention
is paid on the case of a vector bundle. In vector bundles, a bijective
correspondence between the mentioned general frames and pairs of bases, in the
vector fields over the base and in the sections of the bundle, is proved. The
(3\ndash index) coefficients of a connection in such pairs of frames and their
transformation laws are considered. The covariant derivatives are also
mentioned on that context.

    The theory of normal frames for connections on bundles is considered
in section~\Ref[5]{Sect7}. Subsect.~\Ref[5]{Subsect7.1} deals with the general case.
 Loosely said, an adapted frame is called normal if the 2\ndash index
coefficients of a connection vanish in it (on some set). It happens that a
frame is normal if and only if it coincides with the frame it is adapted to.
The set of these frames is completely described in the most general case. The
problems of existence, uniqueness, \etc of normal frames adapted to holonomic
frames, \ie adapted to local coordinates, are discussed in
Subsect.~\Ref[5]{Subsect7.2}. If such frames exist, their general form is
described. The existence of frames normal at a given point and/or along an
injective horizontal path is proved. The flatness of a connection on an open
set is pointed as a necessary condition of existence of (locally) holonomic
frames normal on that set. Some links between the general theory of normal
frames and the one of normal frames in vector bundles, presented in
chapter~\ref{Chapt4}, are given in Subsect.~\Ref[5]{Subsect7.3}. It is proved
that a frame is normal on a vector bundle with linear connection if and only if
in it vanish the 3\ndash index coefficients of the connection. The equivalence
of the both theories on vector bundles is established.

	In section~\Ref[5]{SectA} is formulated and proved a necessary and
sufficient condition for existence of coordinates normal along injective
mappings with non\ndash vanishing horizontal component, in particular along
injective horizontal mappings.

	Section~\Ref[5]{Sect8} is devoted to some aspects of the axiomatical
approach to parallel transport theory%
~\cite{Lumiste-1964,Lumiste-1966,Teleman,Dombrowski,Lumiste-1971,
Mathenedia-4,Durhuus&Leinaas, Khudaverdian&Schwarz,Ostianu&et_al.,Nikolov,
Poor} and its relations to connection theory; it is based on the
paper~\cite{bp-C-TP}. It starts with a definition of a transport along paths
in a bundle and a result stating that, under some assumptions, it defines a
connection. The most important properties of the parallel transports generated
by connections are used to be (axiomatically) defined the concept of a parallel
transport (irrespectively to some connection on a bundle). In a series of
results are constructed bijective mappings between the sets of transports along
paths satisfying some additional conditions, connections, and parallel
transports. In this way, two different, but equivalent, systems of axioms
defining the concept ``parallel transport'' will be established.

	The chapter ends with some remarks and conclusions in
Sect.~\Ref[5]{Conclusion}.


\include{nf-book5}

\include{nf-book9}

\end{document}

%% file: nf-book0.tex
\pagenumbering{roman}

\renewcommand{\thefootnote}{\fnsymbol{footnote}} 
\maketitle				
\renewcommand{\thefootnote}{\arabic{footnote}}   

\setcounter{page}{4}		



\typeout{}
\typeout{????????????????????????????????????????????????????????????????}
\typeout{}
\typeout{This is the front matter of the book }
\typeout{"Handbook of normal frames and coordinates"}
\typeout{by Bozhidar Zakhariev Iliev.}
\typeout{Its initial draft version was written during the period}
\typeout{March, 1999 -- ???, 2000}
\typeout{}
\typeout{vvvvvvvvvvvvvvvvvvvvvvvvvvvvvvvvvvvvvvvvvvvvvvvvvvvvvvvvvvvvvvvvv}
\typeout{}


\renewcommand{\contentsname}{Table of Contents}
\tableofcontents		
\listoftables			
\listoffigures	       		

\chapter*{List of conventions}
\markboth{\slshape LIST OF CONVENTIONS}{\slshape LIST OF CONVENTIONS}
\addcontentsline{toc}{chapter}{List of conventions}
\thispagestyle{myheadings}
\markboth{\slshape LIST  OF CONVENTIONS}{\slshape LIST  OF CONVENTIONS}

\paragraph*{References.}

The book is divided into chapters which have a sequential Roman
enumeration. The chapters are divided into sections with a sequential Arabic
enumeration, which is independent in each chapter. Some sections are divided
into subsection,

	In each chapter the subsections, equations, propositions, theorems,
lemmas, and so on have a double independent enumeration of the
form~\textbf{m.n} or~\textbf{(m.n)} for the equations,
\textbf{m,n=1,2,$\boldsymbol{\ldots}$}\ where~\textbf{m} is the number of the
section in which the designated item appears and~\textbf{n} is its sequential
number in it. So, proposition~4.7 and~(3.12) (or equation~(3.12)) mean
respectively proposition~7 in section~4 and equation~12 in section~3 of the
current chapter. A suitable item from a chapter different from the current
one is referred as~\textbf{R.m}, or~\textbf{R.m.n} or~\textbf{(R.m.n)} for
equations, where \textbf{R=I,II,$\boldsymbol{\ldots}$}\ is the Roman number
of the chapter in which the item appears; \eg remark~II.5.3 and IV.4 (or
section~IV.4) mean respectively remark~3 in section~5 of chapter~II and
section~4 in chapter~IV.

	The footnotes are indicated as superscripts in the main text and have
independent Arabic enumeration in each section. When we refer to a footnote,
it is on the current page if the page on which it appears is not explicitly
indicated.

\paragraph*{Citations.}
An Arabic number in square brackets, e.g.~[27], directs the reader to the
list of references, \ie in this example~[27] means the~27$^\mathrm{th}$ item
from the Bibliography list beginning on page~\pageref{Bibliography}.

\paragraph*{The ends of the proofs} are marked by empty square sign,
viz.\ with \QEDsymbol.

\paragraph*{Indices.}	\label{indices}
The Latin indices refer to an arbitrary linear (vector) space, in
particular to the tangent and cotangent spaces. If in a given problem are
presented the tangent and cotangent spaces to a manifold and other vector
space(s), then the indices referring to the first two spaces are denoted with
small Greek letters; for the rest one(s) the Latin letters will be used.

\paragraph*{Einstein's summation convention:}	\Label{Einstein}
in a product of quantities or in a single expression, a summation over indices
repeated on different levels is assumed over the whole range in which they
change. Any exception of this rule is explicitly stated.

\paragraph*{Symmetrization and antisymmetrization.}
	\Label{(anti)symmetrization}%
On indices included in (or surrounded by) round (resp.\ square) brackets a
symmetrization (resp.\ antisymmetrization) with coefficient one over the
factorial of their number is assumed. If some indices in such a group have to
be excluded from this operation, they are included in (surrounded by) vertical
bars.

\paragraph*{Matrix of linear mapping}
with respect to a given basis, or bases, or field of, possibly local,
bases:  the same symbol but the kernel letter is in
\textbf{boldface}. Exception: the matrix of a derivation (derivative
operator) is denoted by boldface capital Greek letter gamma, \ie by
$\bs{\Gamma}$, possibly with some indices.

\paragraph*{Matrix elements.}	\label{MatrixElements}
When the elements of a (two-dimensional) matrix
are labeled by superscript and subscript, the superscript is considered as a
first index, numbering the matrix's rows, and the subscript as a second one,
numbering the matrix's columns. In this way the matrix of composition of
linear mappings is equal to the product of the matrices of the mappings in
the same order in which they appear in the composition and this does not
depend on the way the matrix's indices are situated.

\paragraph*{Free arguments.}	\label{FreeArguments}
If we want to show explicitly the argument(s) of some mapping or to single
out it (them) as arbitrary while the other arguments, if any, are
considered as fixed ones, we denote it (them) by (centered) dot, \ie by
$\cdot$. E.g., if $f\colon A\to C$ and $g\colon A\times B\to C$, then
$f(\,\cdot\,)\equiv f$, $g(\,\cdot\,,\,\cdot\,)\equiv g$, and $g(\,\cdot\,,b)$,
$b\in B$, means $g(\,\cdot\,,b)\colon A\to C$ with
$g(\,\cdot\,,b)\colon a \mapsto g(a,b)$ for all $a\in A$.

\chapter*{Preface}
\markboth{\slshape PREFACE}{\slshape PREFACE}
	\Label{Preface}
	\addcontentsline{toc}{chapter}{Preface}
\thispagestyle{myheadings}
\markboth{\slshape PREFACE}{\slshape PREFACE}

	The main subject of this book is an up-to-date and in-depth survey of
the theory of normal frames and coordinates in differential geometry. The
existing results, as well as new ones obtained lately by the author, on the
theme are presented.

	The text is so organized that it can serve equally well as a
reference manual, introduction to and review of the current research on the
topic. Correspondingly, the possible audience ranges from graduate and
post\ndash graduate students to scientists working in differential geometry and
theoretical/mathematical physics. This is reflected in the bibliography which
consists mainly of standard (text)books and journal articles.

	The present monograph is the first attempt for collecting the known
facts concerting normal frames and coordinates into a single publication. For
that reason, the considerations and most of the proofs are given in details.

	Conventionally local coordinates or frames, which can be
holonomic or not, are called \emph{normal} if in them the coefficients of a
linear connection vanish on some subset, usually a submanifold, of a
differentiable manifold. Until recently the existence of normal frames was
known (proved) only for symmetric linear connections on submanifolds of a
manifold. Now the problems concerning normal frames for derivations of
the tensor algebra over a differentiable manifold are well investigate; in
particular they completely cover the exploration of normal frames for
arbitrary linear connections on a manifold. These rigorous results are
important in connection with some physical applications. They may be applied
for rigorous analysis of the equivalence principle. This results in two
general conclusions: the (strong) equivalence principle (in its
`conventional' formulations) is a provable theorem and the normal frames are
the mathematical realization of the physical concept of `inertial' frames.
The normal frames find other important physical application in the bundle
formulation of quantum mechanics. It turns out that in a normal frame the
bundle Heisenberg and Shcr{\"o}dinger pictures of motion. coincide.

	Applying some freedom of language, we can state the general  physical
idea: the normal frames are the most suitable ones for describing free
objects and events, i.e.\ such that on them do not act any forces. Regardless
of the different realizations of that idea in general relativity and its
generalizations, quantum mechanics, gauge theories etc., there is an
underlying mathematical background for the general description of such
situations: the existence (or non\ndash existence) of normal frames in vector
bundles. This observation fixes to a great extend the mathematical tools
required for the description of some fundamental physical theories.

	In the book, formally, may be distinguished three parts: The first one
includes chapters~\Ref{Chapt1}\Ndash\Ref{Chapt3} and deals with a variety of
mathematical problems concerning normal frames and coordinates on
differentiable manifolds. The second part consists of chapters~\Ref{Chapt4}
and~\Ref{Chapt5} and investigates normal frames (and possibly coordinates) in
vector bundles and differentiable bundles, respectively. The last part,
involving the text after chapter~\Ref{Chapt5}, contains inquiry material.

	The requisite mathematical language required for the description of
normal frames is spread over the initial sections of the chapters. In
particular,
sections~\Ref[1]{Sect2}\Ndash~\Ref[1]{Sect4},
\Ref[3]{Sect2}, \Ref[4]{Sect7}, \Ref[4]{Sect8}, \Ref[4]{131Subsect1}
and~\Ref[5]{Sect2}--\Ref[5]{Sect6}
can be collected into an introductory chapter under the title ``Mathematical
preliminaries''~%
\footnote{~%
As (practically) any `preliminary' knowledge requires for its understanding
some other `preliminary' to it knowledge, in the corresponding sections are
cited a number of works containing this second kind of mathematical
`luggage'.%
}
but this is not done by pedagogical reasons.~%
\footnote{~%
The material is so organized, that the required concepts and results appear
in the logical order in which they are necessary for some particular
purpose(s).%
}
	The normal coordinates and frames, in the case of linear connections
on a manifold, are initially introduced in chapter~\Ref{Chapt1}. It
contains our basic preliminary material and a review of the Riemannian
coordinates. Chapter~\ref{Chapt2} is devoted to the existence, uniqueness,
construction and other related problems concerning normal frames and
coordinates in manifolds endowed with linear connection. It presents, in
historical order, a detailed review of the existing literature as well as
generalization of a number of results, \eg for connections with torsion.
	Further, in chapter~\Ref{Chapt3}, problems connected with the
existence, uniqueness, holonomicity \etc of normal frames for arbitrary
derivations of the tensor algebra over a manifold are investigated.
	Next (chapter~\Ref{Chapt4}), the same range of problems is explored
for normal frames for linear transports in vector bundles. This material
covers completely the special case of normal frames for linear connections in
vector bundles or on a differentiable manifold.
	The main aim of chapter~\Ref{Chapt5} is the exploration of normal
frames (and coordinates, if any) for general connections on differentiable
fibre bundles which, in particular, can be vector ones.

	The general approach of the book is essentially coordinate-dependent
or basis\ndash dependent. This is due to its basic subject: frames, bases or
coordinates with some special properties. However, if possible and suitable,
the coordinate\nobreak-free notation and methods are not neglected.

	The basic mathematical prerequisites vary form chapter to chapter but
generically they include the grounds of vector (linear) spaces, differentiable
manifolds, vector bundles, connection theory, and a firm belief in the
existence and uniqueness theorems of ordinary differential equations. Some of
the corresponding concepts and results are reproduced in our text but the
acquaintance with adequate literature is required. Appropriate references
are given in the Introductions to the chapters and directly in the main text.

	The material is so organized that a successive chapter generalizes
the preceding one(s) and refers to it (them).

\begin{center}
\tiny$\star$ \footnotesize$\star$ \large$\star$ \LARGE$\star$
\Huge$\star$ $\bigstar$ $\star$
\LARGE$\star$ \large$\star$ \footnotesize$\star$ \tiny$\star$
\end{center}
\nopagebreak

	Any suggestions and comments are welcome. The author's postal
address is
	\begin{quote}\sffamily\bfseries
Bozhidar Zakhariev Iliev,
Laboratory of Mathematical Modeling in Physics,
Institute for Nuclear Research and Nuclear Energy,
Bulgarian Academy of Sciences,
Boul.\ Tzarigradsko chauss\'ee~72, 1784 Sofia, Bulgaria,
	\end{quote}
his e-mail address is
	\begin{quote}\sffamily\bfseries
bozho@inrne.bas.bg
	\end{quote}
and
	\begin{quote}\sffamily\bfseries
http://theo.inrne.bas.bg/$\sim$bozho/
	\end{quote}
is the URL address of his personal World Wide Web site.

 \vspace*{3.2ex}
\noindent
 \hfill%
\begin{minipage}{8em}
Sofia, Bulgaria\\
7 July, 2006
\end{minipage}

\chapter*{Acknowledgments}
\markboth{\slshape ACKNOWLEDGMENTS}{\slshape ACKNOWLEDGMENTS}
	\Label{Acknowledgments}
	\addcontentsline{toc}{chapter}{Acknowledgments}
\thispagestyle{myheadings}
\markboth{}{\slshape ACKNOWLEDGMENTS}

 	The idea for writing a book on normal frames and coordinates in
differential geometry was suggested by~professor\ Dr.~Michiel~Hazewinkel
(Centre for Mathematics and computer science, Amsterdam, The Netherlands) to
the author in October~1998. I express to him my sincere gratitude for this.

	I would like to thank Prof.\ Dr.\ Stancho Dimiev and
mathematicianVladimir Aleksandrov (Mathematical Institute of Bulgarian Academy
of Sciences) for the encouragement and fruitful conversations at the early
stages of my investigations on normal frames. Special thanks to professor
Dr.~Sava Manov (1943\Ndash 2005) (Institute for Nuclear Research and Nuclear
Energy of Bulgarian Academy of Sciences). During 1991\Ndash 1996 he read the
manuscripts of some my papers on the subject. With him I have had many useful
discussions concerning, first of all, the applications of normal frames to
gravity physics. I am indebted to prof.~N.A.~Chernikov for his hospitality
during my numerous visits in the Joint Institute for Nuclear Research, Dubna,
Russia.

	The (anonymous for me) referees from Journal of Physics~A:
Mathematical and General helped me with many suggestions and improvements. I
would like to thank them also.

	The entire manuscript was typeset by the author by means of the main
\LaTeXe\ document preparation computer system together with a number of
additional to it packages of programs, including first of all AMS-\LaTeX,
AMS-fonts, \textsc{Bib}\TeX, MakeIndex, Index, etc. My gratitude to all of
the numerous persons who created and developed, maintain and distribute
(free!) this valuable and high\ndash quality typesetting system.

	I would like to express my deep gratitude to my mother Dr.~Lilyana
Stefanova Shtereva for her understanding, support, love and faith in me.

	The work on this monograph was partially supported by the National
Science Fund of Bulgaria under Grant No.~F~1515/2005.

\vspace{4ex}

\hfill \emph{Bozhidar Z. Iliev}

\noindent%
11 July, 2006
\\ Sofia, Bulgaria

\newpage\hspace*{0pt}


%
%

\newpage
\pagenumbering{arabic}
\setcounter{page}{1}

%% file: nf-book9.tex
\typeout{}
\typeout{????????????????????????????????????????????????????????????????}
\typeout{}
\typeout{These are the reference and index lists of the book }
\typeout{"Handbook of normal frames and coordinates"}
\typeout{by Bozhidar Zakhariev Iliev.}
\typeout{Their initial draft version was written during the period}
\typeout{March, 1999 -- July, 2000}
\typeout{}
\typeout{vvvvvvvvvvvvvvvvvvvvvvvvvvvvvvvvvvvvvvvvvvvvvvvvvvvvvvvvvvvvvvvvv}
\typeout{}

\bibliography{bozhopub,bozhoref}
\bibliographystyle{unsrt}


 	\index{local!coordinate system|see{coordinate system}}
 	\index{local!coordinates|see{coordinates}}
 	\index{differentiable!mapping|see{mapping differentiable}}
\index{tangent! field|see{tangent vector field}}
\index{vector field|see{\emph{also} tangent vector field}}
\index{tangent! mapping|see{differential of mapping}}
\index{one-form|see{covector field}}
\index{1-form|see{covector field}}
\index{covariant vector field|see{covector field}}
 	\index{curvature free |see{flat}}
 	\index{torsion free |see{torsionless}}
 	\index{flat!linear connection|see{\emph{also} linear connection}}
 	\index{torsionless!linear connection|see{\emph{also} linear connection}}
	\index{transport!parallel|see{parallel transport}}
\index{derivative!covariant|see{covariant derivative}}
\index{derivative!Lie|see{Lie derivative}}
\index{derivative!Fermi-Walker|see{Fermi-Walker derivative}}
\index{derivative!Fermi|see{Fermi derivative}}
\index{derivative!Truesdell|see{Truesdell derivative}}
\index{derivative!Jaumann|see{Jaumann derivative}}
\index{differentiation|see{derivation}}
\index{parallel tensor field|see{tensor field}}
\index{geodesics|see{geodesic path}}
\index{autoparallels|see{autoparallel path}}
\index{equation of!geodesics|see{geodesic equation}}
\index{mapping!exponential|see{exponential mapping}}
\index{manifold!Riemannian|see{Riemannian}}
	\index{frame!normal|see{normal frame, normal frame for}}
	\index{normal frame|see{\emph{also} normal frame for}}
	\index{coordinates!normal|see{normal coordinates, normal coordinates
for and \emph{also} normal frame}}
	\index{normal coordinates|see{\emph{also} normal coordinates for, normal frame,
normal frame \emph{holonomic},}}
	\index{normal coordinates!on \ldots|see{normal coordinates on \ldots}}
	\index{chart!normal|see{normal chart}}
	\index{normal chart|see{\emph{also} normal coordinates}}
\index{components' matrix of|see{matrix of}}
\index{coefficients' matrix of|see{matrix of}}
\index{linear connection!of mixed type|see{mixed linear connection}}
	\index{normal coordinates for linear connection!at a
point|see{\emph{also} Riemannian normal coordinates}}
	\index{normal coordinates for linear connection!along
path|see{\emph{also} Fermi coordinates}}
\index{transport!along paths!linear|see{linear transport along paths}}
\index{derivation along paths!normal frame for|see{normal frame for}}
\index{derivation along tangent vector fields!normal frame for|see{normal frame for}}
 \index{total bundle space|see{bundle space}}
 \index{standard fibre|see{fibre of bundle}}
 \index{typical fibre|see{fibre of bundle}}
 \index{derivation along paths in vector bundle|see{derivation along paths}}
 \index{vector bundle!derivation along paths in|see{derivation along paths}}
 \index{bundle!tangent to manifold|see{tangent bundle}}
 \index{bundle!tensor|see{tensor bundle}}
 \index{linear connection!derivation along paths generated by|see{derivation
along paths}}
	\index{covariant derivative operator|see{covariant derivative}}
\index{section-curvature operator of|see{curvature of}}
\index{curvature operator of|see{curvature of}}
\index{torsion vector|see{torsion}}
\index{holonomic normal frame|see{normal frame \emph{holonomic}}}
	\index{parallel transport!in tangent bundle!along paths|see{parallel
transport along paths}}
\index{normal frame!strong|see{strong normal frame}}
	\index{system of parallel transport|see{parallelism structure}}
	\index{parallel lifting |see{lifting}}
	\index{parallel section|see{section}}
	\index{parallel path|see{path}}
	\index{vector bundle!connection on|see{connection}}
\index{curvature tensor field|see{curvature}}
\index{Euclidean!linear transport|see{linear transport}}
\index{Euclidean!derivation|see{derivation}}
\index{equation of!normal frames|see{normal frame equation}}
\index{connection!linear|see{linear connection}}
\index{covariant derivative!in vector bundle|see{\emph{also} linear connection}}
\index{coordinates normal|see{\emph{also} normal coordinates}}
\index{manifold!Minkowski|see{Minkowski spacetime}}
\index{Levi-Civita connection|see{\emph{also} Riemannian connection}}
\index{equation of!normal coordinates|see{normal coordinates equation}}
\index{vector|see{\emph{also} tangent vector}}



\newpage\vspace*{0pt}

\printindex[authors]

\newpage\vspace*{0pt}

\printindex[symbols]

\printindex[subject]
